\documentclass[a4paper]{article}

\usepackage{delarray,verbatim,enumerate}
\usepackage{amsmath,amsthm,amstext,amsbsy,amssymb,amsfonts,amscd}
\usepackage[all]{xy}
\usepackage[frenchb]{babel}

\usepackage{aeguill}

\newtheorem{Th}{Th\'eor\`eme}[]
\newtheorem{Lem}[Th]{Lemme}

\def\Remarque{\smallskip\noindent {\it Remarque.~}}

\font\teneufm=eufm10
\font\seveneufm=eufm7
\font\fiveeufm=eufm5
\newfam\gothfam
\textfont\gothfam=\teneufm \scriptfont\gothfam=\seveneufm
\scriptscriptfont\gothfam=\fiveeufm
\def\goth{\fam\gothfam}

\def\F2{\mathbb{F}_2}		\def\Q{\mathbb Q}		\def\Z{\mathbb Z}	
 	\def\z2{\mathbb{Z}_2} 	

 				\def\U{\mathcal  U}
\def\J{\mathcal  J}  	\def\C{\mathcal  C}	\def\R{\mathcal  R}
 	  	\def\Cl{\mathcal  C \ell}
	 	\def\Cs{\mathcal C{\ell}s}

		\def\p{{\goth p}}		

\def\wi{\widetilde}		
\def\sg{\operatorname{sg}}	\def\deg{\operatorname{deg}}
\def\Gal{\operatorname{Gal}}	\def\Log{\operatorname{Log}}

\begin{document}

\title{\bf\Large Sur le sous-groupe des éléments de hauteur infinie du
$K_2$ d'un corps de nombres\footnote{Acta Arithmetica {\bf 122} (2006), 235--244.}}

\author{Jean-François {\sc Jaulent} \& Florence {\sc Soriano-Gafiuk}}
\date{}
\maketitle

{\footnotesize
\noindent{\bf Résumé.} Poursuivant l'approche logarithmique des noyaux de la $K$-théorie pour les corps de nombres (cf. \cite{J1, J2, JS1, JS2}), nous déterminons ici en toute généralité les 2-rangs respectifs du noyau sauvage $W\!K_2(F)$ et du sous-groupe $K_2^\infty(F)$ des éléments de hauteur infinie dans $K_2(F)$, à l'aide des groupes de classes positives du corps $F$.\par

Le résultat principal de cette étude, qui conjugue les critères d'égalité de $W\!K_2(F)$  et de $K_2^\infty(F)$ obtenus par K. Hutchinson (cf. \cite{H1, H2, HR}) avec la description logarithmique des noyaux sauvages, permet ainsi de déterminer effectivement en termes de classes logarithmiques dans quels cas le sous-groupe $K_2^\infty(F)$  est un facteur direct du noyau sauvage $W\!K_2(F)$.\par}
\medskip

{\footnotesize
\noindent{\bf Abstract.} By using the logarithmic approach of the classical kernels for the $K_2$ of number fields (cf. \cite{J1, J2, JS1, JS2}), we compute the 2-rank of the wild kernel $W\!K_2(F)$ and the 2-rank of the subgroup $K_2^\infty(F) = \cap_{n\ge 1}\, K_2^n(F)$ of infinite heigh elements in $K_2(F)$ in terms of positive class groups for any  number field $F$. \par 
The main result, which connects the criterion of equality between these two groups already obtained by K. Hutchinson (cf. \cite{H1, H2, HR}) and the logarithmic description of these groups, gives a simple logarithmic characterization of the number fields $F$ for which the subgroup $K_2^\infty(F)$ is a direct factor of the wild kernel $W\!K_2(F)$.\par}

\


\bigskip

\noindent{\bf Introduction}

\medskip
  On sait depuis les travaux de Tate que, pour tout corps de nombres $F$, le noyau $W\!K_2(F)$ 
des symboles sauvages donnés par la théorie locale du corps de classes est
un sous-groupe fini de $K_2(F)$ dont l'arithmétique est mystérieusement reliée à celle
des groupes de classes d'idèles ou d'idéaux (cf. e.g. \cite{Ta}, \cite{ Ko}). \smallskip

C'est ainsi que J. Tate a montré que pour tout premier impair $\ell$, la $\ell$-partie de $W\!K_2(F)$ coïncide avec celle du sous-groupe $K_2^\infty(F) = \cap_{n \in \mathbb N^\times}\, K_2^n(F)$ des éléments de hauteur infinie dans $K_2(F)$~; et il est montré dans \cite{JS1} que le $\ell$-rang commun de ces deux groupes est donné par l'isomorphisme de modules galoisiens~:\smallskip

\centerline{${}^{\ell}K_2^\infty(F) = {}^{\ell}W\!K_2(F) \simeq {}^\Delta(\mu_\ell \otimes \wi\Cl_{F[\zeta_\ell]})
= \mu_\ell \otimes \wi\Cl_{F[\zeta_\ell]}^{e_{\bar\omega}}$,}\smallskip

\noindent où ${}^{\ell}W\!K_2(F) = W\!K_2(F) /W\!K_2(F)^\ell$ est le quotient d'exposant $\ell$ de $W\!K_2(F)$~; $\mu_\ell$ désigne le groupe des racines $\ell$-ièmes de l'unité~; $F[\mu_\ell]$ est l'extension cyclotomique de $F$ correspondante~; $\Delta = \Gal(F[\mu_\ell]/F)$ le groupe de Galois associé~; $\wi\Cl_{F[\zeta_\ell]}$ est le $\ell$-groupe des classes logarithmiques du corps $F[\mu_\ell]$~; et ${}^\Delta(\mu_\ell \otimes \wi\Cl_{F[\zeta_\ell]})$ est le groupe des copoints fixes du tensorisé de $\wi\Cl_{F[\zeta_\ell]}$ par $\mu_\ell$, qui s'identifie encore, comme expliqué dans \cite{JS1}, au tordu par $\mu_\ell$ de la composante anticyclotomique du groupe $\wi\Cl_{F[\zeta_\ell]}$. C'est la situation algorithmiquement exploitée dans \cite{C1} et généralisée aux noyaux étales supérieurs dans \cite{JM}.\medskip

Lorsque $\ell$ est égal à 2, la situation est plus complexe à double titre~:\pagebreak

- d'abord parce que les 2-parties des deux groupes $W\!K_2(F)$ et $K_2^\infty(F)$ ne coïncident plus nécessairement, le sous-groupe des éléments de hauteur infinie dans $K_2(F)$ pouvant être contenu stictement dans le noyau sauvage $W\!K_2(F)$~;

- ensuite parce que les isomorphismes logarithmiques précédents peuvent être en défaut dès lors que le 
 corps $F$ ne contient pas les racines 4-ièmes de l'unité comme l'attestent les calculs menés dans \cite{So2} pour certains corps quadratiques.\medskip

De fait, l'examen plus attentif de la formule explicite pour le symbole sauvage $\Big( \frac{-1,\ x}{F_{\p}} \Big)$ établie dans \cite{J1} montre qu'interviennent dans ce cas non seulement les valuations logarithmiques attachées aux places finies, mais aussi les fonctions signes attachées aux places réelles ou 2-adiques, ce qui justifie l'introduction de la notion de classes signées (cf. \cite{So1}, \cite{J4}), qu'on peut regarder comme l'analogue logarithmique de celle classique de classes au sens restreint.\smallskip

Les résultats de \cite{JS2} permettent ainsi de décrire la 2-partie du noyau sauvage à l'aide d'un quotient convenable du groupe des classes signées dit {\em groupe des classes positives}. L'isomorphisme obtenu \smallskip

\centerline{${}^2W\!K_2(F) \simeq \{\pm 1\} \otimes \Cl^{\,pos}_F$}\smallskip

\noindent laisse cependant ouverte la question de la description du groupe $K_2^\infty(F)$ et de son quotient ${}^2 K_2^\infty(F)$. L'objet de la présente note est précisément de combler cette lacune en décrivant ce dernier groupe en termes de classes positives.
\smallskip

\bigskip

\noindent{\bf 1. Rappels sur les places logarithmiquement signées ou primitives}

\bigskip

Commençons par rappeler quelques définitions indispensables à la compréhension du Théorème principal~:\medskip

\noindent {\bf Définition} (cf. \cite{H1}, \cite{JS2}, \cite{Ry}, \cite{H3}, \cite{HR}). Soit $F$ un corps de nombres~; 
pour chaque place $\p$ de $F$, écrivons $F_\p$ le complété de $F$ en $\p$. Cela posé~:
\medskip

\noindent \quad {\it (1)} Nous disons d'une part que la place $\p$ est~:\smallskip

{\it \ (i) }{\it signée}, lorsque le corps $F_\p$ ne contient pas les racines 4-ièmes 
de l'unité~;\smallskip

{\it (ii) logarithmiquement signée} (ou encore {\it exceptionnelle}), lorsque la $\z2$-exten\-sion cyclotomique $F_\p^c$ de $F_\p$ ne contient pas les racines 4-ièmes de l'unité\footnote{De sorte que ne sont {\it logarithmiquement signées} que les places réelles et certaines places paires, i.e. au-dessus de 2.}~;\smallskip

{\it (iii) logarithmiquement primitive}, lorsqu'elle ne se décompose pas dans (le premier
étage de) la $\z2$-extension cyclotomique $F^c$ de $F$.\medskip

\noindent \quad {\it (2)} Nous disons d'autre part que le corps $F$ est~:\smallskip

{\it \ (i) }{\it signé}, lorsqu'il ne contient pas les racines 4-ièmes de l'unité (ce qui a lieu si 
et seulement s'il possède au moins une place {\it signée})~;\smallskip

{\it (ii) exceptionnel}, lorsque sa $\z2$-extension cyclotomique $F^c$ ne 
contient pas les racines 4-ièmes de l'unité~;\smallskip

{\it (iii) logarithmiquement signé}, lorsque sa 2-extension abélienne localement 
cyclotomique maximale $F^{lc}$ ne contient pas les racines 4-ièmes de l'unité 
(ce qui a lieu si et seulement s'il possède au moins une place {\it logarithmiquement 
signée})~;\smallskip

{\it (iv) logarithmiquement primitif } lorsque $F$ possède au moins une place (nécessairement paire) qui est à la fois {\it logarithmiquement signée et primitive}~; \medskip

\Remarque Il résulte des définitions précédentes qu'un corps de nombres est {\em logarithmiquement signé} lorsqu'il est {\em localement exceptionnel} en l'une de ses places~; il est alors (globalement) {\it exceptionnel}, mais la réciproque est fausse (cf. infra). En revanche un corps de nombres est {\em signé} si et seulement s'il l'est en l'une de ses places (auquel cas il l'est en une infinité d'entre elles en vertu du théorème de \v Cebotare\v v).\medskip

Donnons un exemple~:\medskip

\noindent {\bf Exemple}. Soit $F=\Q[\sqrt d]$ (avec $d \in \mathbb Z , \ d \ne 0$ ou 1, $d$ sans facteur carré) un corps quadratique réel ou imaginaire. On a alors les équivalences~:\smallskip

 {\it \ (i) } $F$ {\it signé} (i.e. $i \notin F$) $\Leftrightarrow$ $d \ne -1$.\smallskip

{\it (ii)} $\ F$ {\it exceptionnel} (i.e. $i \notin F^c$) $\Leftrightarrow$ 
$d \ne -1,-2$.\smallskip

{\it (iii)} $F$ {\it logarithmiquement signé} (i.e. $i \notin F^{lc}$) 
$\Leftrightarrow$ ($d> 0$, i.e. $F$ est {\it $\infty$-signé}) ou ($d \not \equiv -1$ 
[mod 8] et $d \not \equiv -2$ [mod 16], i.e. $F$ est 2-{\it logarithmiquement signé} ).\smallskip

{\it (iv)} $\ F$ {\it logarithmiquement primitif}\, $\Leftrightarrow$ ($d \not\equiv -1$ 
[mod 8] et $d \not\equiv \pm 2$ [mod 16]) ou $d=2$ (les places paires sont à la fois {\it logarithmiquement signées et primitives}).\medskip

Introduisons maintenant quelques notations~: pour chaque place non complexe $\p$ du corps considéré $F$, désignons par \smallskip

\centerline{$\R_{F_\p} = \varprojlim F_\p^\times/F_\p^{\times 2^n}$}\smallskip

\noindent le compactifié 2-adique du groupe multiplicatif du complété de $F$ en $\p$, que la théorie locale du corps de classes identifie au groupe de Galois de la pro-2-extension abélienne maximale de $F_\p$.\smallskip 

L'approche logarithmique amène à considérer deux applications définies sur le compactifié 
2-adique $\R _{F_\p}$ et à valeurs respectivement dans $\Z_2$ et dans $\{ \pm
 1\}$~:\smallskip

$\qquad \bullet$ la $\p$-{\it valuation logarithmique} $\wi v_\p$ d'une part, donnée par la formule~:\smallskip

\centerline{$\wi v_\p (\cdot )\, =\, -\frac{\Log\mid \cdot \mid _{\p}}{\deg \p}$\quad (cf. \cite{J2}, p. 306, d\'ef. 1.3 ($iv$)),}\smallskip

\noindent où $\mid \cdot \mid _{\p}$ désigne la valeur absolue 2-adique à valeurs dans le groupe multiplicatif $\Z_2^\times$ et $\deg \p$ est un facteur de normalisation choisi pour qu'on ait $\wi v_\p(\R_{F_\p})=\Z_2$~;\smallskip

$\qquad \bullet$ la {\it fonction signe } $\sg_\p$ d'autre part, dont une expression est~:\smallskip

\centerline{$\sg_\p (\cdot )\, =\,  \varepsilon (\mid \cdot \mid _{\p})$\quad (cf. \cite{J4}, p. 457, déf. 1),}\smallskip

\noindent où $\varepsilon$ désigne la projection naturelle de $\Z_2^\times\, \simeq\, \{\pm 1\}\times (1+4\Z_2)$ sur $\{\pm1\}$. \smallskip

Le noyau $\wi \U_{F_\p}$ de la valuation $\wi v_\p$ est, par définition, le sous-groupe des {\it unités logarithmiques} de $\R_{F_\p}$~: c'est le groupe de normes associé par la Théorie locale du corps de classes à la $\Z_2$-extension cyclotomique $F_\p^c$ de $F_\p$ (avec la convention $F_\p^c = F_\p$ aux places réelles)~; celui de la fonction signe $\sg_\p$ est le sous-groupe des {\it éléments positifs} $\R^+_{F_\p}$~: il correspond à l'extension $F_\p[i]$ engendrée par les racines 4-ièmes de l'unité~; l'intersection $\wi \U^+_{F_\p} = \wi \U_{F_\p} \cap \R^+_{F_\p}$, enfin,  est le sous-groupe des {\it unités logarithmiques positives}~: et il correspond à l'extension composée $F^c_\p[i]$, i.e. à la 2-extension cyclotomique de $F_\p$.\smallskip

En résumé, la place $\p$ est donc {\it signée} lorsque la fonction signe $\sg_\p$
est non triviale, auquel cas le sous-groupe positif $\R^+_{F_\p}$ est d'indice 2 
dans $\R_{F_\p}$~; elle est {\it logarithmiquement signée} lorsque la fonction
signe $\sg_\p$ est non triviale sur le sous-groupe unité $\wi \U_{F_\p}$, auquel cas
$\wi \U^+_{F_\p}$ est d'indice 2 dans $\wi \U_{F_\p}$. Et cette dernière éventualité ne se 
produit qu'aux places réelles et à certaines des places au-dessus de 2. Nous disons que ce sont les places paires logarithmiquement signées.

\bigskip

\noindent{\bf 2. Compléments sur les groupes de classes positives}

\medskip

Introduisons maintenant les groupes de classes positives et rappelons d'abord pour cela quelques notations de la Théorie du corps de classes~:\smallskip

Etant donné un corps de nombres $F$, désignons par $F^{ab}$ sa pro-2-extension {\it abélienne maximale}~; notons $F^c$ la {\em $\Z_2$-extension cyclotomique} de $F$ et $F^{lc}$ la sous-extension {\em localement cyclotomique} maximale de $F^{ab}$, i.e. la plus grande sous-extension de $F^{ab}$ qui est complètement décomposée sur $F^c$\smallskip

Par la Théorie 2-adique du corps de classes (cf. \cite{J3}), le groupe de Galois $\Gal(F^{ab}/F)$ s'identifie comme groupe topologique au quotient $\C_F=\J_F/\R_F$ du 2-adifié du groupe des idèles du corps $F$, défini comme le produit restreint\smallskip

\centerline{$\J_F = \prod_\p^{res} \R_{{F_\p}}$ ,}\smallskip

\noindent par son sous-groupe principal $\R_F = \Z_2 \otimes_\Z F^\times$. En d'autres termes, dans la correspondance du corps de classes, le groupe $\J_F$ est associé à $F$ et le sous-groupe $\R_F$ à $F^{ab}$. Il se trouve que le groupe de normes attaché à la $\Z_2$-extension cyclotomique $F^c$

\centerline{$\wi\J_F = \{(x_\p)_\p \in \J_F\mid \, \sum_\p \wi v_\p(x_\p) \deg \p = 0 \}$}\smallskip

\noindent diffère (en général) de celui attaché à l'extension localement cyclotomique $F^{lc}$\smallskip

\centerline{$\wi\U_F \R_F = \prod_\p \wi\U_{F_\p} \R_F$ ,}\smallskip

\noindent de sorte que le 2-groupe des classes logarithmiques (de degré nul) \smallskip

\centerline{$\wi \Cl_F = \wi\J_F /\wi\U_F \R_F \simeq \Gal(F^{lc}/F^c)$ ,}\smallskip

\noindent mesure l'écart entre les 2-extensions {\it localement} et {\it globalement} cyclotomiques $F^{lc}$ et $F^c$. On ne sait pas encore si ces derniers groupes sont finis pour tout les corps de nombres $F$ (ce que postule précisément une généralisation naturelle de la conjecture de Gross), mais ce point est sans importance pour les calculs numériques, puisque les groupes de classes logarithmiques se calculent essentiellement comme les groupes de classes ordinaires et que la validité de la conjecture, qui revient à affirmer la non trivialité d'un certain régulateur, se vérifie aisément en pratique dès que l'on sait calculer dans le corps étudié (cf. \cite{C1}).\smallskip

Lorsque le corps $F$ est {\em exceptionnel}, il faut en outre distinguer la $\Z_2$-extension cyclotomique $F^c$ de la 2-tour cyclotomique $F^c[i]$. Le groupe de normes attaché à $F[i]$ est le noyau de la formule du produit pour les fonctions signes~:\smallskip

\centerline{$\J^*_F = \{(x_\p)_\p \in \J_F\mid \, \prod_\p \sg_\p(x_\p) = +1 \}$~;}\smallskip

\noindent et coïncide avec le produit $\J^+_F \R_F = \prod_\p \R^+_{{F_\p}} \R_F$, d'après la caractérisation locale de l'extension quadratique $E = F[i]$. L'intersection\smallskip

\centerline{$\wi\J_F^* = \wi\J_F \cap \J_F^* = \{(x_\p)_\p \in \J_F\mid \, \prod_\p \mid x_\p \mid_\p = 1 \}$}\smallskip

\noindent correspond donc au compositum $E^c = F^c[i]$. Et la plus grande extension localement triviale de $E^c$ qui est abélienne sur $F$, en d'autres termes l'intersection $F^{ab} \cap E^{lc}$, est ainsi associée au sous-groupe de normes~:\smallskip

\centerline{$\wi\U_F^+ \R_F = \prod_\p \wi\U^+_{{F_\p}} \R_F = N_{E/F}(\wi\U_E) \R_F$.}\smallskip

\noindent Le quotient correspondant (cf. \cite{J4}) \smallskip

\centerline{$\wi{\Cl s}_F = \wi\J_F^*/ \wi\U_F^+ \R_F \simeq \Gal (E^{lc} \cap F^{ab} /E^c)$}\smallskip

\noindent est, par définition, le 2-{\em groupe des classes logarithmiques signées} du corps $F$~; c'est aussi l'image, par la norme arithmétique $N_{E/F}$ du 2-groupe des classes logarithmiques de $E$, comme on le voit sur le schéma de corps~:

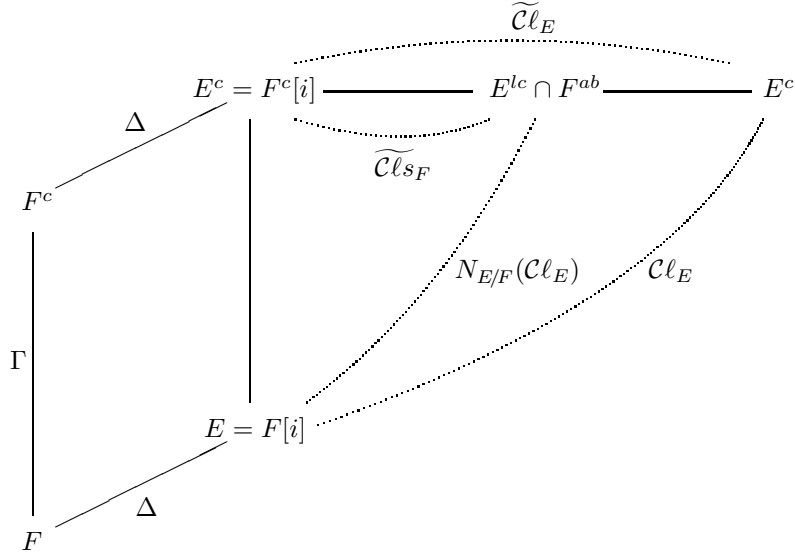
\begin{center}

\unitlength=1.5cm
\begin{picture}(7.2,5)

\put(0.1,0){$F$}
\put(0.2,0.3){\line(0,1){2.5}}
\put(0.1,3){$F^c$}
\put(0,1.6){$\Gamma$}

\put(1.7,1){$E = F[i]$}
\put(2.1,1.3){\line(0,1){2.5}}
\put(1.6,4){$E^c=F^c[i]$}

\put(0.4,0.2){\line(2,1){1.5}}
\put(1.1,0.3){$\Delta$}
\put(0.4,3.2){\line(2,1){1.5}}
\put(1,3.7){$\Delta$}

\put(2.75,4.05){\line(1,0){1.3}}
\put(4.2,4){$E^{lc} \cap F^{ab}$}
\put(5.2,4.05){\line(1,0){1.3}}
\put(6.6,4){$E^c$}

\bezier{80}(2.5,4.3)(4.4,4.7)(6.3,4.3)
\put(4.4,4.6){$\wi\Cl_E$}

\bezier{40}(2.5,3.8)(3.5,3.5)(4.2,3.8)
\put(3.2,3.3){$\wi{\Cl s}_F$}

\bezier{80}(2.6,1.3)(3.8,2.2)(4.6,3.8)
\put(3.9,2.4){$N_{E\!/\!F}(\Cl_E)$}

\bezier{120}(2.7,1.1)(5.8,2.2)(6.6,3.8)
\put(5.6,2.4){$\Cl_E$}

\end{picture}
\end{center}\medskip

Il est commode de noter \smallskip

\centerline{$\Cl s_F = \J^*_F/ \wi\U^+_F \R_F = N_{E/F}(\Cl_E)$}\smallskip

\noindent le groupe analogue pris sans condition de degré, image par la norme arithmétique du groupe des classes logarithmiques (toujours sans condition de degré) $\Cl_E$.
Avec ces notations, le 2{\em -groupe des classes positives} $\Cl^{\,pos}_F$ est le quotient \smallskip

\centerline{$\Cl^{\,pos}_F = \J^*_F / \J^{pos}_F \R_F$} \smallskip

\noindent du groupe $\Cl s_F$ par l'image du sous-groupe $\J^{pos}_F = \prod_{\p\in PLS}\R_{{F_\p}}^+\prod_{\p\notin PLS}\wi\U_{{F_\p}}^+$ de $\J^*_F$ construit à partir de l'ensemble $PLS$ des places de $F$ logarithmiquement signées. Et le sous-groupe des {\em classes positives de degré nul} est le quotient\smallskip

\centerline{$\wi{\Cl}^{\,pos}_F = \wi\J^*_F / \wi\J^{pos}_F \R_F$,} \smallskip

\noindent du groupe $\wi{\Cl s}_F$, image dans $\Cl^{\,pos}_F$ du sous-groupe $\wi\J^*_F$ de $\J^*_F$.\smallskip

Via le degré, le quotient $\,\Cl^{\,pos}_F / \wi\Cl{}^{\,pos}_F \simeq \J^*_F/\wi\J^*_F\prod_{\p \in PLS}\R^+_{F_\p}$ s'identifie à $\Z_2$ en l'absence de places paires logarithmiquement signées, auquel cas le sous-groupe de degré nul est toujours un facteur direct de $\,\Cl^{\,pos}_F$. En revanche, en présence de places paires logarithmiquement signées, c'est un quotient fini de $\Z_2$, et dans ce cas il peut arriver que le sous-groupe de degré nul soit ou ne soit pas facteur direct dans $\,\Cl^{\,pos}_F$ ~: s'il l'est, le quotient ${}^2\wi\Cl{}^{\,pos}_F$ est alors un hyperplan du $\F2$-espace vectoriel ${}^2\Cl{}^{\,pos}_F$~; s'il ne l'est pas, on a~:  ${}^2\wi\Cl{}^{\,pos}_F \simeq {}^2\Cl{}^{\,pos}_F$.

\bigskip

\noindent{\bf 3. \'Enoncé du Théorème principal}

\medskip

Ces notations étant données, nous pouvons maintenant énoncer le résultat principal de cette note, qui précise ceux obtenus par J.-F. Jaulent \cite{J1} sur le $\ell$-noyau sauvage, puis par K. Hutchinson \cite{H1}, \cite{H2} ou J.-F. Jaulent \& F. Soriano \cite {JS2} sur le 2-noyau sauvage~:

\begin{Th} Soient $F$ un corps de nombres, $W\!K_2(F)$ le noyau des symboles sauvages et 
$K_2^\infty(F) = \cap_{n \ge 1} K_2^n(F)$ le sous-groupe des éléments de hauteur infinie dans $K_2(F)$ (qui est d'indice au plus 2 dans $W\!K_2(F)$). Cela posé~:\smallskip

(i)\, Si le corps $F$ n'est pas exceptionnel (i.e. pour $i \in F^c$), on a 
directement~:\smallskip

\centerline{$K_2^\infty(F) / K_2^\infty(F)^2 = W\!K_2(F)/W\!K_2(F)^2 \ 
\simeq \ \{ \pm 1 \} \otimes _{\z2} \wi \Cl_F$,}\smallskip 

\noindent où $\,\wi \Cl_F$ désigne le 2-groupe des classes logarithmiques 
(au sens ordinaire) de $F$.

\par Et, dans ce cas, on a toujours l'égalité~: $W\!K_2(F) = K_2^\infty(F)$.\smallskip

(ii) Si le corps $F$ est exceptionnel (i.e. pour $i \not\in F^c$), on a canoniquement~:
\medskip

\centerline{$W\!K_2(F)/W\!K_2(F)^2 \ \simeq \ \{ \pm 1 \} \otimes _{\z2} 
\Cl^{\,pos}_F$,}\smallskip

\noindent où  $\,\Cl^{\,pos}_F$ désigne le 2-groupe des classes 
logarithmiques positives de $F$, mais~:\smallskip

\centerline{$K_2^\infty(F) / K_2^\infty(F)^2 \ \simeq \ \{ \pm 1 \} 
\otimes _{\z2} \wi\Cl{}^{\,pos}_F$,}

\noindent où $\,\wi\Cl{}^{\,pos}_F$ désigne le sous-groupe de 
$\,\Cl{}^{\,pos}_F$ formé des classes de degré nul.

\par Et, dans ce cas, on a : $W\!K_2(F) = K_2^\infty(F)$ si et 
seulement si $F$ est logarithmiquement primitif~; et $(W\!K_2(F) : K_2^\infty(F)) = 2$ sinon.
Enfin, dans cette toute dernière situation, le sous-groupe $K_2^\infty(F)$ est facteur direct du groupe $W\!K_2(F)$ si et seulement si le sous-groupe des classes de degré nul $\,\wi\Cl{}^{\,pos}_F$ est lui même facteur direct du groupe $\,\Cl^{\,pos}_F$ des classes positives.
En d'autres termes, si $F$ est exceptionnel et logarithmiquement imprimitif, on a~:\par
- soit \quad $\,\,\,\Cl^{\,pos}_F \simeq \,\,\wi\Cl{}^{\,pos}_F \oplus \Z/2^s\Z$ \quad et \quad $W\!K_2(F) \simeq K_2^\infty(F)) \oplus \F2$~;\par
- soit \quad ${}^2\Cl^{\,pos}_F \simeq {}^2\wi\Cl{}^{\,pos}_F$ \quad et \quad ${}^2W\!K_2(F) \simeq {}^2K_2^\infty(F))\, (\simeq {}^2\Cl^{\,pos}_F)$.
\end{Th}

\Remarque Lorsque le corps $F$ est {\it exceptionnel}, mais non {\it logarithmiquement
signé}, i.e. lorsque l'on a $i \in F^{lc} \setminus F^c$, le groupe des idèles positifs $\J^{pos}_F$ défini dans \cite{JS2} coïncide avec le groupe des unités logarithmiques $\wi \U_F$ et le quotient $\J_F / \J^{pos}_F
\R_F = \J_F / \wi \U_F \R_F = \Cl_F$ est alors le 2-groupe des classes logarithmiques (sans condition de degré) du corps $F$. En particulier, les quotients d'exposant 2 respectifs 
${}^2 \Cl^{\,pos}_F$ et ${}^2 \wi \Cl_F$ du groupe des classes positives $\Cl^{\,pos}_F = 
\J^*_F / \J^{pos}_F \R_F$ et du groupe des classes logarithmiques (de degré nul) $\wi \Cl_F =
\wi \J_F / \wi \U_F \R_F$ sont les hyperplans noyaux dans le $\F2$-espace vectoriel
${}^2 \Cl_F$ des formes linéaires induites respectivement par le signe $\sg_F$ et le degré $\deg_F$. Ils sont donc (non canoniquement) isomorphes, et l'isomorphisme annoncé $^2 W\!K_2(F) \simeq
{}^2 \Cl^{\,pos}_F$ redonne bien dans ce cas l'isomorphisme $^2 W\!K_2(F) \simeq {}^2 \wi 
\Cl_F$ établi dans \cite{JS2}.

\bigskip

\noindent{\bf 4. Preuve du Théorème principal}

\medskip

Le cas où le corps considéré $F$ n'est pas {\it exceptionnel} ne pose pas problème~:  nous avons alors $W\!K_2(F) = K_2^\infty (F)$ d'après \cite{H1}~; et $^2 W\!K_2(F) \simeq {}^2 \wi \Cl_F$ d'après \cite{J1} si le corps est {\it signé}, d'après \cite{JS2} s'il ne l'est pas~; d'où les isomorphismes annoncés dans ce cas.\smallskip

Supposons donc le corps $F$ {\it exceptionnel}, introduisons l'extension $E = F[i]$ et considérons le schéma d'extensions dessiné plus haut. Dans cette dernière situation, nous avons, en vertu du Lemme 3.2 de \cite{H2}~:

\begin{Lem} L'homomorphisme de transfert $Tr_{E\!/\!F}$ de $K_2(E)$ dans $K_2(F)$ envoie le noyau sauvage $W\!K_2(E) = K_2^\infty(E)$ sur le sous-groupe $K_2^\infty(F)$ d'indice au plus 2 dans 
$W\!K_2(F)$ formé des éléments de hauteur infinie dans $K_2(F)$.
\end{Lem}

Et simultanément, comme expliqué dans la section 2~:

\begin{Lem} L'opérateur norme $N_{E\!/\!F}$ envoie le 2-groupe des classes logarithmiques $\,\wi \Cl_E$ sur le groupe des classes signées $\,\wi \Cs_F$~; donc, par passage au quotient, le groupe ${}^2 \wi\Cl_E$ {\it sur} l'image canonique du groupe $\,{}^2 \wi\Cl {}^{\,pos}_F$ dans $\,{}^2 \Cl^{\,pos}_F$.
\end{Lem}

Considérons donc le diagramme commutatif ci-dessous qui relie les isomorphismes entre les quotients d'exposant 2 respectifs des groupes de classes et des noyaux sauvages écrits en haut pour le corps $E$ et en bas pour le corps $F$~:

\begin{displaymath}\xymatrix{
{}^2 \wi\Cl_E \ = \ {}^2 \wi \Cl{}^{\,pos}_E \ \ar@{->>}[d]_{N_{E\!/\!F}} \ar@{->>}[r]^{\sim \phantom{20}}_{\beta_E \phantom{20}} & \ {}^2 W\!K_2(E) = {}^2 K^\infty_2(E) \ar@{->>}[d]^{Tr_{E\!/\!F}} \\
{}^2 \wi{\Cl s}_F \ = \ {}^2 N_{E\!/\!F}(\wi \Cl_E) \ \ar@{->>}[d]_\pi \ar@{->>}[r]_{\delta_F \phantom{20}} & \ {}^2 K^\infty_2(F) \ar@{->>}[d]^\gamma \\
{}^2 \wi\Cl{}^{\,pos}_F \ar[d]_\alpha \ar@{->>}[r]_{\eta_F \phantom{20}}  & K^\infty_2(F) W\!K_2(F)^2 / W\!K_2(F)^2 \ar[d]^\iota \\
{}^2 \Cl^{\,pos}_F \ar@{->>}[r]_{\beta_F \phantom{20}}^{\sim \phantom{20}}  & {}^2 W\!K_2(F)
}\end{displaymath}

Distinguons les deux cas~:\medskip

$\ (i)$ Si le corps $F$ est logarithmiquement primitif, le groupe des classes positives $\,\Cl^{\,pos}_F$ coïncide avec son sous-groupe de degré nul $\,\wi\Cl{}^{\,pos}_F$ comme expliqué plus haut~;  ainsi $\alpha$ est l'identité et $\iota$ est surjectif. En particulier, il suit $W\!K_2(F) = K^\infty_2(F)$ (conformément au résultat de \cite{H1})~;  et nous obtenons, comme annoncé dans ce cas, les isomorphismes naturels~:\smallskip

\centerline{$K_2^\infty(F) / K_2^\infty(F)^2 = W\!K_2(F)/W\!K_2(F)^2 \ 
\simeq \ \{ \pm 1 \} \otimes _{\Z_2} \Cl_F^{\, pos}$,}\medskip

$(ii)$ Si le corps $F$ n'est pas logarithmiquement primitif, le sous-groupe des classes de degré nul $\,\wi\Cl{}^{\,pos}_F$ est strictement contenu dans le groupe des classes positives $\,\Cl^{\,pos}_F$ et il peut en être ou non facteur direct~:\smallskip

Si $\,\wi\Cl{}^{\,pos}_F$ est facteur direct de $\,\Cl^{\,pos}_F$, l'application $\alpha$ est injective et son image est un hyperplan du $\F2$-espace vectoriel ${}^2 \Cl^{\, pos}$. Il suit que $\eta_F$ est elle-même injective (donc bijective) et que le sous-groupe $K^\infty_2(F) W\!K_2(F)^2 / W\!K_2(F)^2$ est un hyperplan du quotient ${}^2 W\!K_2(F)$. En particulier $K^\infty_2(F)$ est lui-même un facteur direct non trivial de $W\!K_2(F)$, et nous avons, comme attendu~:\smallskip

\centerline{${}^2 K_2^\infty(F) \simeq \ \{ \pm 1 \} \otimes _{\Z_2}\wi\Cl{}_F^{\, pos} \quad \& \quad  {}^2 W\!K_2(F) \simeq \ \{ \pm 1 \} \otimes _{\Z_2} \Cl_F^{\, pos} $.}\medskip

Inversement, si le sous-groupe $K^\infty_2(F)$ est un facteur direct non trivial de $W\!K_2(F)$, l'application $\iota$ n'est pas surjective et l'application $\alpha$ ne l'est pas non plus~; il en résulte que $\,\wi\Cl{}^{\,pos}_F$ est alors un facteur direct de $\,\Cl^{\,pos}_F$.\smallskip

En fin de compte, les morphismes $\alpha$ et $\iota$ sont simultanément surjectifs ou non~; les sous-groupes $\,\wi\Cl{}^{\,pos}_F$ et $K^\infty_2(F)$ sont simultanément facteurs directs ou non respectivement de $\,\Cl^{\,pos}_F$.\ et de $W\!K_2(F)$. Et, bien entendu, lorsqu'ils ne le sont pas, nous avons les isomorphismes non canoniques~:\smallskip

\centerline{${}^2 K_2^\infty(F) \simeq \ {}^2 W\!K_2(F) \simeq \ \{ \pm 1 \} \otimes _{\Z_2} \Cl_F^{\, pos} \simeq \ \{ \pm 1 \} \otimes _{\Z_2}\wi\Cl{}_F^{\, pos}$.}

\Remarque Il suit de là que, contrairement à ce qui est indiqué malheureusement dans \cite{JS2}, l'application {\it canonique} ${}^2\wi\Cl_F^{\, pos} \rightarrow {}^2\Cl{}_F^{\, pos}$ n'est pas, alors, bijective.\medskip

Diverses illustrations numériques sont données dans un travail en cours \cite{C2}, qui généralise dans ce contexte signé l'approche algorithmique initiée dans \cite{C1}.


{\small
\def\refname
}
\bigskip\noindent

\begin{tabular}{l p{0.8 cm}  l}
Jean-François {\sc Jaulent}	&{}	& Florence {\sc Soriano-Gafiuk}\\
Institut de Mathématiques 	&{} 	& Département de Mathématiques\\
Université  Bordeaux 	1	&{} 	& Université de Metz\\
 351, cours de la Libération	&{} 	& Ile du Saulcy\\
 F-33405 TALENCE  Cedex	&{}	& F-57045 METZ Cedex\\
jaulent@math.u-bordeaux1.fr 	&{} 	& soriano@poncelet.univ-metz.fr
\end{tabular}

\end{document}